\documentclass[preprint,12pt,authoryear,comma,times]{elsarticle}
\usepackage{amsmath,amssymb,graphicx,pifont}

\usepackage{lineno}
\usepackage[usenames]{color}
\definecolor{lightgray}{gray}{0.2}
\newcommand{\degr}{$^{\circ}$}  

\biboptions{sort&compress}
\journal{Process Safety and Environmental Protection}
\begin{document}

\begin{frontmatter}


\title{Oscillatory thermal instability and the Bhopal disaster}

\author{R. Ball}

\address{Mathematical Sciences Institute, The Australian National University\\ Canberrra ACT 0200 Australia }
\ead{Rowena.Ball@anu.edu.au}
\ead[url]{www.maths.anu.edu.au/$\sim$ball/}

\begin{abstract}
A stability analysis is presented of the hydrolysis of methyl isocyanate (MIC) using a homogeneous flow reactor paradigm.  The results simulate the thermal runaway that occurred inside the storage tank of MIC at the Bhopal Union Carbide plant in December 1984. The stability properties of the model indicate that the thermal runaway may have been due to a large amplitude, hard thermal oscillation initiated at a subcritical Hopf bifurcation. This type of thermal misbehavior cannot be predicted using conventional thermal diagrams, and may be typical of liquid thermoreactive systems.  
\end{abstract}

\begin{keyword}
Methyl isocyanate \sep Thermal runaway \sep Oscillatory instability \sep Bhopal \sep Liquid explosives

\end{keyword}

\end{frontmatter}
\linenumbers
\section{Introduction \label{sec1}}
More than 25 years after the Bhopal disaster its horrific legacy is now well-documented \citep{Mishra:2009}, but the causes are still being debated in the international media. Was the tragedy due to neglect, parsimony, or procrastination by Union Carbide on safety and maintenance? Ignorance, corruption, sabotage and cover-up? Inadequate regulation of urban and industrial development? Possibly all of the above, but they are putative, secondary or socioeconomic contributing factors. (A brief account of the disaster is given in the Appendix.) The primary cause of the thermal runaway that led to the venting of a poisonous mist of methyl isocyanate over the city of Bhopal was physicochemical.  In this belated work I present a simple stability analysis of the thermokinetics of methyl isocyanate (MIC) hydrolysis, revealing rogue thermal misbehavior that appears to be endemic to reactive organic liquids and that cannot be predicted using conventional heat generation/loss rate diagrams. 

The Bhopal incident is by no means the only example of disastrous thermal runaway occurring in a storage tank. The Seveso accident in Italy in July 1976 in which large quantities of a toxic dioxin were released into the environment occurred under storage tank conditions~\citep{Theophanous:1983} , while more or less minor runaways and explosions due to unforseen reactions in storage tanks and vessels are relatively common.  The special dangers of chemical storage were discussed by~\cite{Gygax:1988}, who pointed out that circumstances favouring heat accumulation are actually \textit{more} likely to occur in storage tanks and equipment parts that are not actively controlled, than in dedicated, well-designed chemical reactors. 

Despite the acknowledged dangers of large-scale storage of reactive chemicals very little has been published on the dynamics of processes that may occur in storage tanks. \cite{Velo:1996} summarize the literature on theoretical and experimental validations of runaway criteria and parametric sensitivity  in batch reactors and storage tanks. In defining critical conditions they, along with  other authors cited therein, begin with the assumption that storage tanks of small volume can be modelled as well-stirred batch reactors with linear thermal coupling to the environment. 

However batch reactors have no non-trivial steady states, and there is no general theory for determing whether a thermal excursion will grow or decay. Here it is assumed that the same parameters that govern the stability of a thermoreactive process in nonequilibrium steady state also govern the stability of thermoreactive processes in storage tanks. From the comprehensive stability and bifurcation analysis of the CSTR that was carried out in \cite{Ball:1999} these parameters are ambient temperature, residence time, heat loss, and intensive properties of the reacting system. It is  shown in this work that a simple spatially homogeneous, steady state approximation can simulate thermoreactive processes in a storage tank with high fidelity. 

Another driver for better understanding of thermoreactive processes in liquids has emerged recently; this is the use of organic hydroperoxide explosives by terrorists.

\section{ \label{sec2}MIC chemistry and known relevant data} 

A chemical analysis of the residue in the MIC storage tank (Tank 610) sampled seventeen days after the event found a variety of condensation products, mainly the cyclic trimer \citep{Dsilva:1986}.  However, experiments to elucidate the organic chemistry of the formation of these products indicated that these condensations must have been initiated at temperatures and pressures well above the normal boiling point of MIC. Therefore, it is thought that the initial reaction of thermokinetic significance was hydration to the unstable N-methyl carbamate, indicated in grey:\\
\centerline{CH$_3$NCO + H$_2$O $\overset{k(T)}\longrightarrow$ \textcolor{lightgray}{CH$_3$NHCOOH},}\\
where $k(T)$ is the temperature-dependent (pseudo-)first order rate constant.  The primary thermal effect that led to the onset of critical conditions is thought to be due to the overall reaction\\
\centerline{CH$_3$NCO + H$_2$O $\overset{k(T)}\longrightarrow$ {\color{lightgray}CH$_3$NHCOOH}
$\longrightarrow$ CH$_3$NH$_2$ + CO$_2$.}

Some relevant physicochemical data and quantities for MIC hydrolysis in Tank 610 are given in table \ref{table1}. The thermodynamic parameters are taken from standard tables. 

\begin{table}\caption{\label{table1}Known relevant data for the MIC thermal runaway at Bhopal.}
\small
\begin{tabular}{p{0.5\textwidth}p{0.4\textwidth}}
\hline\\
Boiling point of MIC at 1\,atm& 39.1\degr\,C\\
Density of MIC& 0.9599\,g/cm$^3$ at 20\degr\,C\\
Constant pressure heat capacity of MIC& 67.7\,J/(K\,mol) = 1188\,J/(K\,kg) \\
Reaction enthalpy for MIC hydrolysis&$-$65.1\,kJ/mol\\
Reaction frequency for hydration of MIC \citep{Castro:1985}&$3.9\times 10^{12}$s$^{-1}$\\
Activation energy for hydration of MIC \citep{Castro:1985}& 64\,kJ/mol\\
\hline\\
Mass of MIC in Tank 610&41 tonnes\\
Initial temperature inside Tank 610& 13\degr\,C\\
Estimated time to criticality& 4 hours\\
\hline
\end{tabular}
\end{table}

\section{Thermokinetics of MIC hydrolysis and thermal instability\label{sec3}}
The spatially homogeneous flow reactor in which a reactant undergoes a first order, exothermic conversion is a simple but elucidatory model for thermoreactive systems when it is appropriate to ignore convection, because as a dynamical system it has non-trivial steady states that can be analysed for stability. The mass and enthalpy equations are 
\begin{align}
V\frac{dc}{dt}&=-Vck(T)+F(c_f-c)\label{e1}\\
V\overline{C}\frac{dT}{dt}&=(-\Delta H)Vck(T) - (F\overline{C} + L)(T -T_a).\label{e2}
\end{align}
$V$ is the reacting volume, $c$ is the reactant concentration, $c_f$ is the reactant concentration in the inflow, $F$ is the volumetric flow rate, $\overline{C}$ is the averaged volumetric specific heat, $T$ is the reaction temperature, $\Delta H$ is the reaction enthalpy, $L$ is the linear heat transfer coefficient, $T_a$ is the ambient temperature.  The temperature-dependent (pseudo)-first order reaction rate constant is
 \begin{align}
 k(T)&=A\exp(-E/RT), \label{e3}
 \end{align}
where $A$ is the reaction frequency, $E$ is the activation energy, and $R$ is the universal gas constant. For numerical and comparative reasons  it is useful to work with the following dimensionless system corresponding to equations  (\ref{e1}--\ref{e2}), using (\ref{e3}):
\begin{align}
\frac{dx}{d\tau}&=-xe^{-1/u}+f(1-x)\label{e4}\\
\varepsilon\frac{du}{d\tau}&=xe^{-1/u} -(\varepsilon f+\ell)(u-u_a),\label{e5}
\end{align}
where $x\equiv c/c_f$, $\tau\equiv tA$, $u\equiv RT/E$, $f\equiv F/VA$, $\varepsilon\equiv\overline{C}E/c_f(-\Delta H)R$, $\ell\equiv LE/c_fVA(-\Delta H)R$, $u_a\equiv RT_a/E$. Numerical analysis of equations (\ref{e4}--\ref{e5}) was carried out using rate, thermochemical and temperature data from table \ref{table1} and values of the inverse residence time $f$, heat loss $\ell$, and inflow concentration $c_f$ were assigned on the basis of available data. 

\begin{figure}[ht]
\centerline{
\includegraphics[scale=0.7]{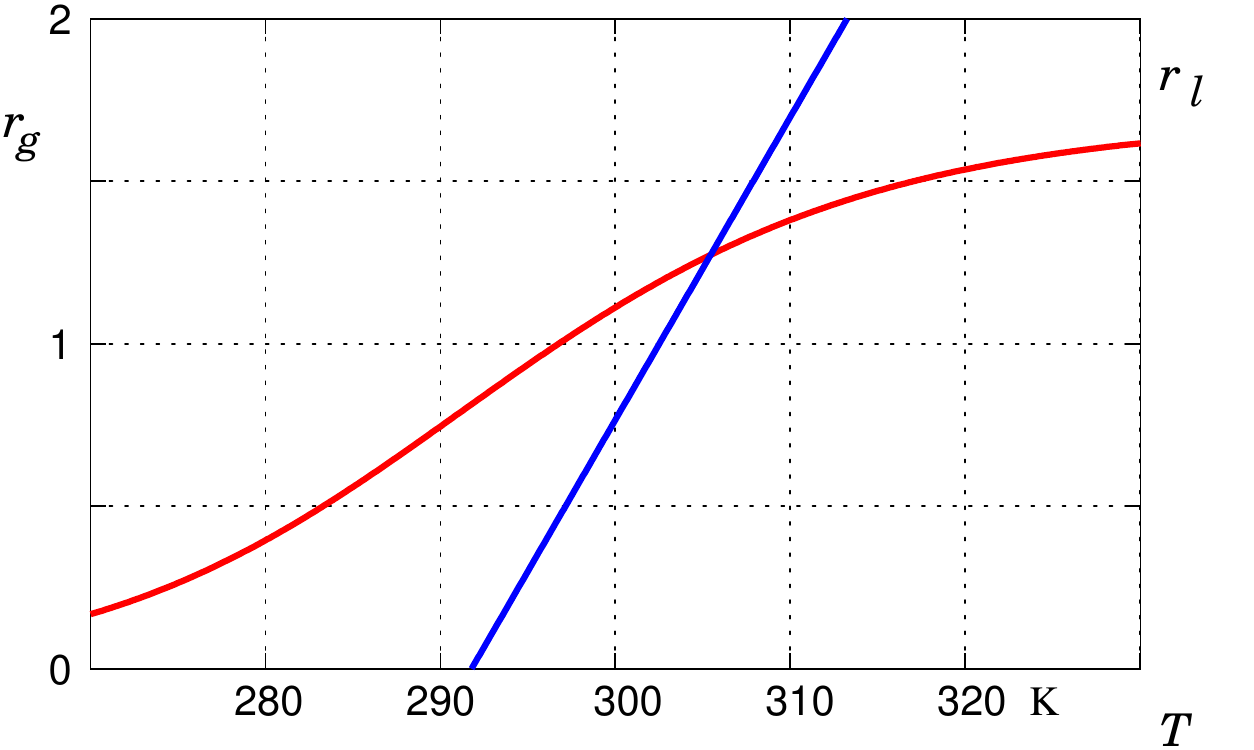}}
  \caption[]{\label{figure1}Rates of reactive heat generation $r_g$ (red) and heat loss $r_l$ (blue) versus  $T$ from equations (\ref{e4}--\ref{e5}). $u_a=0.0379$ (corresponding to $T_a=292\,K $), $f=1.7$, $\ell=700$, $\varepsilon=10$. }
\end{figure}

The reacting mixture self-heats if the rate of reactive heat generation $r_g$ exceeds the linear cooling rate $r_l$. Thermal runaway occurs if $r_g$ exceeds $r_l$ beyond a system-specific threshold; for the hydrolysis of MIC this is taken as the boiling point of MIC. The steady-state rates from equations (\ref{e4}--\ref{e5}), are plotted in figure~\ref{figure1}, where the temperature is labeled in dimensional units. 
The reaction self-heats until the reactor temperature $T$ reaches the steady state temperature of $\sim$305\,K at which the heating and cooling rates are balanced. Since the boiling point of MIC is 312\,K, on the basis of this diagram we would not expect a thermal runaway to develop, even when the ambient temperature is allowed to drift slowly up to 292\,K. 

\begin{figure}[ht]
\centerline{
\includegraphics[scale=0.8]{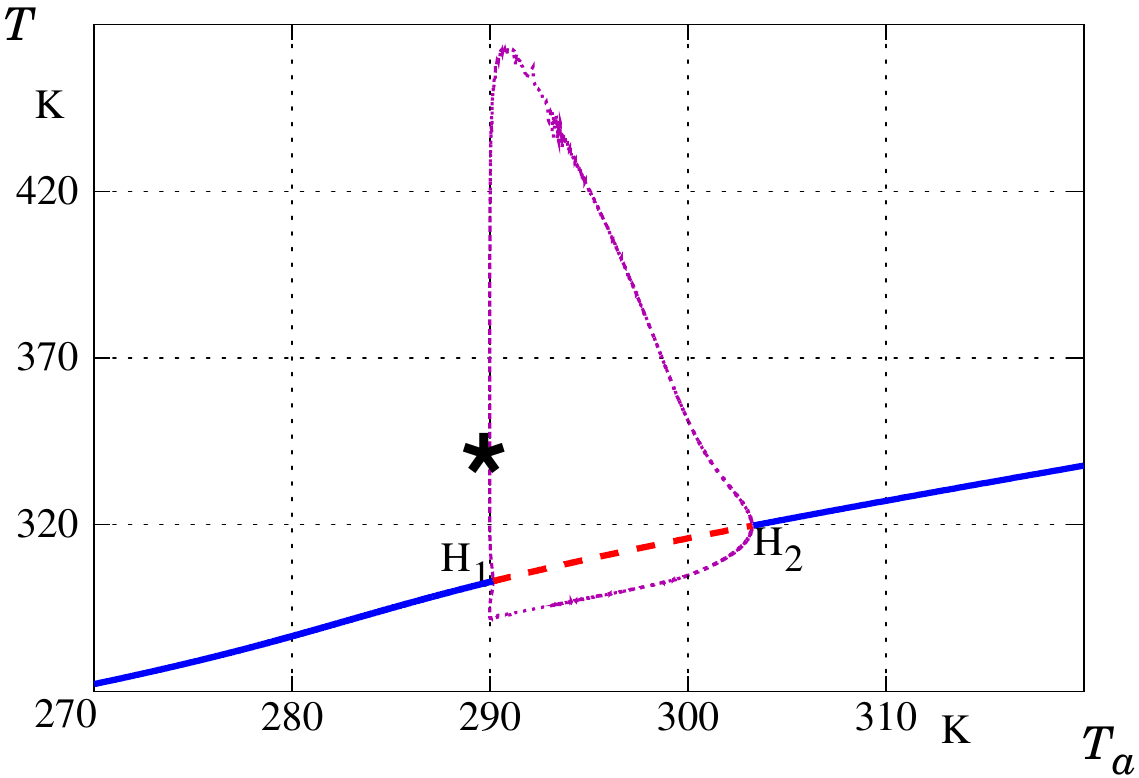} }
  \caption[]{\label{figure2} Bifurcation diagram. Stable steady states are plotted with solid blue line, unstable steady states with dashed red line, and the amplitude envelope of periodic solutions is marked with thin dotted magenta line. $H_1$ and $H_2$ label the Hopf bifurcation and the large  \textbf{\ding{81}} marks the change in stability of the limit cycles. $f=1.7$, $\ell=700$, $\varepsilon=10$.}
\end{figure}

However thermal balance diagrams such as that in figure \ref{figure1} can be dangerously misleading because they infer stability rather than assess stability rigorously, although such diagrams are often used in chemical reactor engineering. The steady states, periodic solutions, and stability analysis of equations (\ref{e4}--\ref{e5}) were computed numerically \citep{Doedel} and yielded a dramatically different picture of the the thermal stability of MIC hydrolysis. Figure \ref{figure2} is a bifurcation diagram in which the steady states and the amplitude envelope of periodic solutions are plotted as a function of $T_a$. 
The steady state is stable at $T_a\approx 286\,K$, the temperature at which the tank of MIC had been held for several months. As $T_a$ is quasistatically increased the reaction temperature $T$ increases slowly, but at $T_a=290.15\,K$ the stability analysis flags an abrupt change in the nature of the solutions. At this point the steady state solutions lose stability to a Hopf bifurcation and the hard onset of a high amplitude thermal oscillation ensues. Clearly, at $T_a=292\,K $ we have catastrophic thermal runaway, contrary to the inadequate prediction given by figure \ref{figure1}. (In the resulting superheated liquid the exothermic condensation reactions would increase the temperature even further.) 

\begin{figure}[ht]
\centerline{
\includegraphics[scale=0.8]{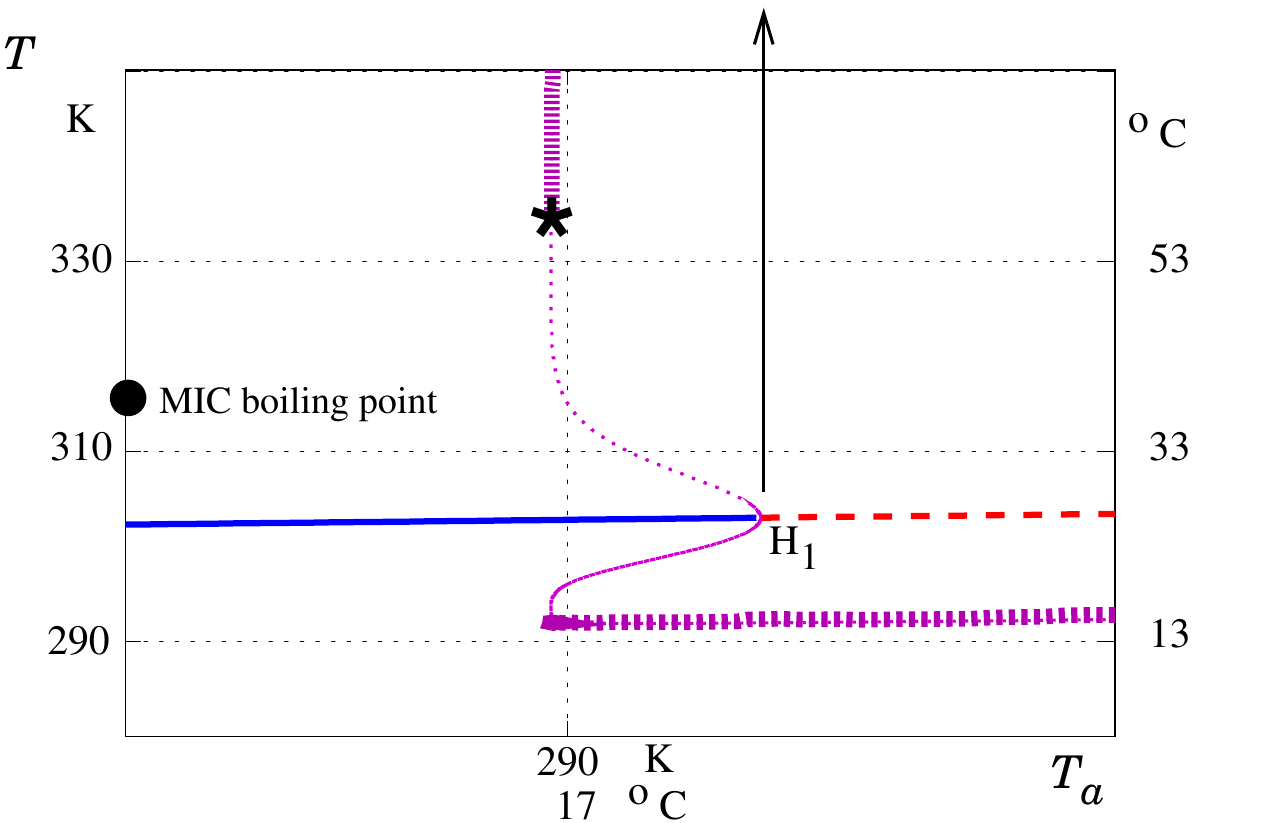}}
  \caption[]{\label{figure3} Close-up of the region around the Hopf bifurcation $H_1$ in figure \ref{figure2}. }
\end{figure}

This is quite different from classical ignition of a thermoreactive system, which occurs at a steady-state turning point. The dynamics of oscillatory thermal runaway can be understood by studying a close-up of the region around the Hopf bifurcation $H_1$ in figure \ref{figure2}. This is shown in figure \ref{figure3}. $H_1$ is subcritical and the emergent limit cycle is \textbf{un}stable. The amplitude envelope of the unstable limit cycles is marked with a thin dotted line; they grow as $T_a$ is \textbf{de}creased. At the turning point of the periodic solution branch marked with a large asterisk the limit cycles become stable. Thermal runaway \textit{may} occur if there are significant perturbations while $T_a$ is within the regime of unstable limit cycles, and it \textit{must} occur when $T_a$ drifts above $H_1$. The arrow indicates the rapid thermal excursion, in principle to the stable limit cycle but in reality the reactant has vaporised, the pressure has soared beyond the safety limits of the tank, and the system must vent since the peak temperature is far above the boiling point of MIC. 

\section{Discussion}
\begin{figure}[ht]
\centerline{
\includegraphics[scale=0.8]{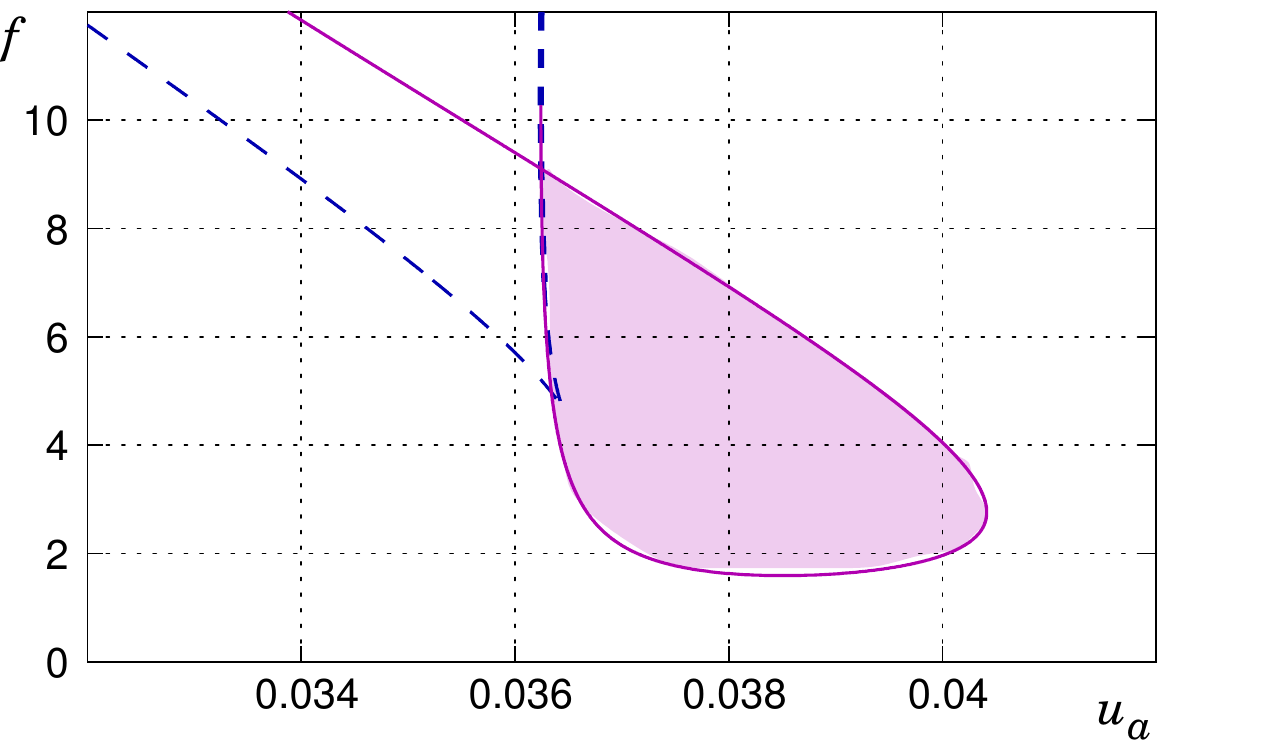}}
  \caption[]{\label{figure4} The locus of Hopf bifurcations is marked with a solid line, the locus of steady-state turning points is marked with dashed line.  $\ell=700$, $\varepsilon=10$. }
  \end{figure}

Although in principle classical thermal ignition is possible for MIC hydrolysis, in practice  the presence of oscillatory instability is all-pervasive and dominant in this system. This can be appreciated by inspection of figure \ref{figure4}, a plot of the loci of the steady state turning points and the Hopf bifurcations of equations (\ref{e4}--\ref{e5}) over the two parameters $u_a$ and the inverse residence time $f$. 
In the filled region thermal runaway will always be oscillatory. The bistable regime occurs at very high flow rates (short residence times). Here classical ignition at a steady state turning point to a stable steady state may occur, but the oscillatory instability is still present and oscillatory thermal runaway occurs (perversely) as the ambient temperature $u_a$ is reduced.  

The tendency for oscillatory, non-classical thermal runaway may be typical of low-boiling, exothermically reactive organic liquids. In the work of \cite{Ball:1995a} the hydration of 2,3-epoxy-1-propanol in a CSTR was found to exhibit similar non-classical thermal misbehavior. The design of safe storage systems for such liquids should focus on damping the oscillatory instability, rather than shifting a classical ignition point using over-simplistic heat generation/loss rate diagrams like figure~\ref{figure1}. 

\begin{figure}[ht]
\centerline{
\includegraphics[scale=0.8]{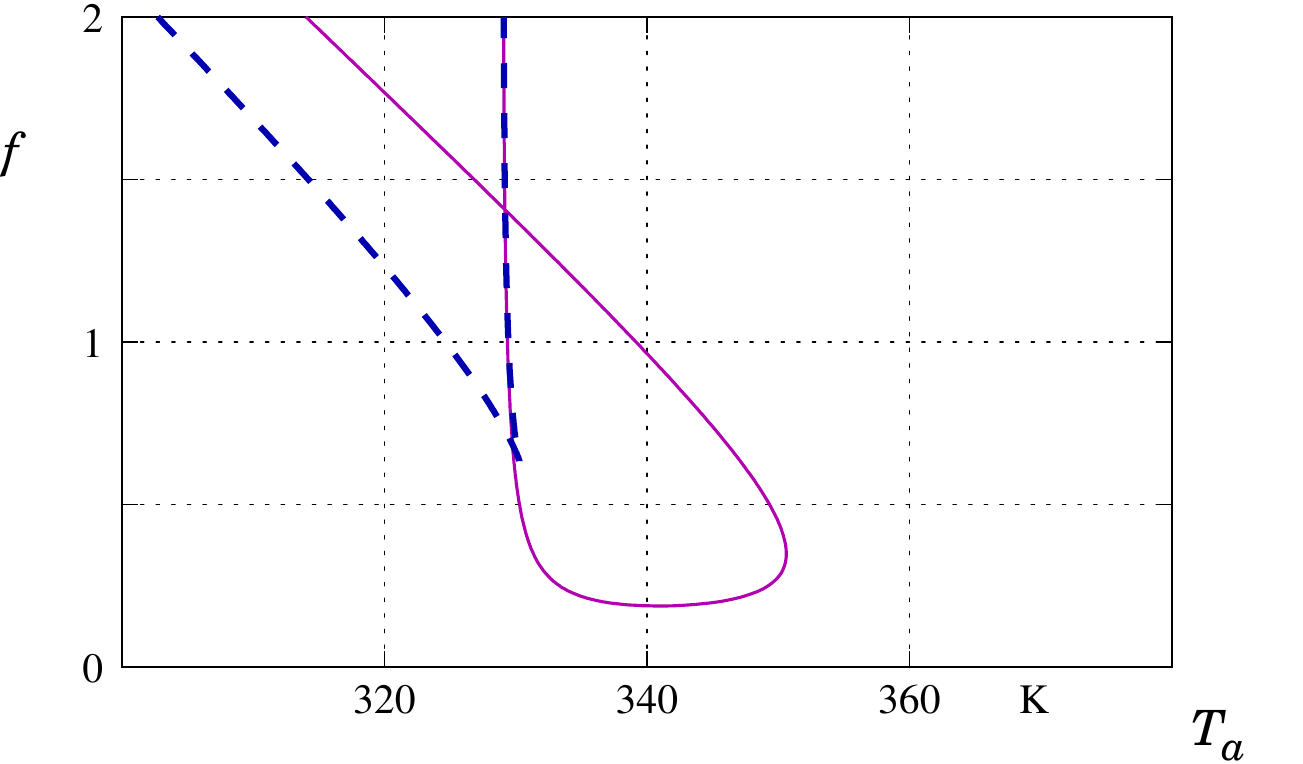}}
  \caption[]{\label{figure5} The 2-parameter bifurcation diagram from equations (\ref{e4}--\ref{e5}) applied to the decomposition of cumene hydroperoxide also has a large region within which thermal runaway is due to oscillatory instability. The locus of Hopf bifurcations is marked with a solid line, the locus of steady-state turning points is marked with dashed line.  $\ell=700$, $\varepsilon=20$. }
  \end{figure}

In recent years low-boiling, exothermically reactive liquids such as a variety of organic peroxides have been used as explosives by terrorists, and their potential use as murder weapons on aircraft is the reason for current restrictions on the liquids that passengers may carry on-board. It is possible that peroxide explosions are due to oscillatory instability rather than classical thermal ignition. Figure \ref{figure5} shows the two-parameter bifurcation diagram using thermokinetic parameters for the decomposition of cumene hydroperoxide \citep{Wu:2010} in equations (\ref{e4}--\ref{e5}). 
In understanding this type of explosion and for deactivating improvised explosive devices that employ such liquids it is clearly necessary to understand oscillatory thermal instability.  

Is it realistic to model a reacting volume inside Tank 610 as a well-stirred flow reactor? Yes, on a timescale over which the reacting volume remains relatively constant and gradientless relative to the much faster rate of reaction. Consider a volume $V$ of liquid in the tank, near the source of the water ingress, within which the hydration reaction takes place.  The reacting volume $V$ grows as water flows into it and reaction proceeds.  This volume can be regarded as the ``reactor'', while the uncontaminated, non-reacting MIC in the rest of the tank can be regarded as the ``coolant'' to which heat is transferred linearly. The reacting volume expands until the rate of inflow matches the rate of outflow of reactant and products from the reacting volume. 

In other words, for the purposes of this model in which the focus is on the dynamics we circumscribe a volume in which the spatial gradients are insignificant in comparison to the time evolution, and therefore can be neglected.  If this approximation does not hold, then we are free to reduce the circumscribed volume until it does. There is nothing particularly artificial or manipulative in doing this; it is just a simplest case scenario for which the powerful tools of stability and bifurcation theory can be applied. 

The Union Carbide plant at Bhopal has been derelict since the disaster and MIC is no longer in use as a bulk chemical, but we cannot afford to close the book on potential problems in storage tanks. In chemical plants and storage sites around the world there are tanks containing reactive organic liquids that have similar physical properties to MIC and undergo reactions with similar thermokinetics. Many other chemicals are stored as bulk commodities that can undergo thermally unstable polymerization, oxidation or decomposition reactions. Analysis of simple gradientless models of processes that may occur in storage tanks provides insight into the physical basis of rogue thermal behavior and a starting point for improved design of safe storage systems. Concomitantly there is a need for more, and more accurate, thermokinetic data to use in such models. It is notable that the data for MIC hydration/hydrolysis in table \ref{table1}  were obtained as part of an unrelated research project initiated before the Bhopal disaster. 

\section{Conclusions\label{sec5}}

\begin{enumerate}
\item A belated investigation of the thermal runaway in Tank 610 that led to the Bhopal disaster has been carried out by modelling the hydrolysis of MIC in the well-stirred limit and analysing the stability of solutions of the dynamical system model. 
\item Thermal runaway occurs due to the onset of a hard thermal oscillation at a subcritical Hopf bifurcation. Classical thermal ignition at a steady state turning point may occur in principle but over the thermal regime of interest the system is dominated by oscillatory instability. 
\item This non-classical oscillatory thermal misbehavior may be generic in reactive organic liquids. 
The results yield valuable intelligence about the causes of themal runaway that may inform improved designs of storage systems for large quantities of toxic and reactive substances. 
\item These results may also inform better management of organic peroxide based explosives. 
\end{enumerate}
\clearpage
\appendix
\section*{Appendix}
The following brief account of the Bhopal disaster has been compiled from the following sources:
\cite{Forman:1985,Weir:1987,Shrivastava:1987,Lepkowski:1994} and \cite{Abbasi:2005}.  

The Union Carbide plant at Bhopal carried out the production of carbaryl, an agricultural insecticide that has been used widely throughout the world since 1945. Methyl isocyanate, a low-boiling, highly reactive and extremely toxic liquid used in the synthesis of carbaryl, was stored in an underground stainless steel tank (Tank 610) which was encased in a concrete shell.  The temperature of the 41 tonnes of MIC in Tank 610 was 12--14\degr\,C rather than the recommended 0--4\degr\,C because the refrigeration unit had been non-operational for several months. On the evening of December 2 1984 a worker had been sent to hose out a nearby tank. The hose was left running unattended, and it is believed that a faulty valve allowed entry of water into the connected Tank 610. (Union Carbide disputes this, asserts that nothing was wrong with its equipment and procedures, and argues that sabotage by a disaffected employee must have caused the disaster.) By 11:30 pm, when workers detected lachrymose whiffs of leaking MIC, water had been running into Tank 610 for at least four hours. Although a slow rise in temperature and pressure in the tank had been noted, the early signs of trouble were not acted upon. Shortly after 11:30 pm the contents of the tank reached thermal criticality and began escaping as vapor from the flare tower.

Downwind of the flare tower lay the crowded suburbs and shantytowns. Most of the fluid in the tank streamed from the tower then drifted low over the city and sank and seeped in deathly mist in lungs in eyes, a period to sleep and swift arrest of retreat. More than 3000 lives were claimed immediately and many tens of thousands through the subsequent days and months and years lost their lives or their health to the poison's effects, and the dead are still being counted. \\

\noindent\textit{Acknowledgement:}
This work is supported by Australian Research Council Future Fellowship FT0991007.


\end{document}